\documentclass[12pt,a4paper,leqno,verbatim]{amsart}
\newcounter{minutes}
\divide\time by 60
\newcounter{hours}
\multiply\time by 60 \addtocounter{minutes}{-\time}
\usepackage{amssymb}
\usepackage{hyperref}
\usepackage{graphicx}
\usepackage{color}

\date{}
\newfont{\cyrilic}{wncyr10 scaled 1000}
\title[On the generalized convexity and concavity]
{On the generalized convexity and concavity}

\author[B. A. Bhayo]{Barkat Ali Bhayo}
\address[B. A. Bhayo]{Department of Mathematical Information Technology, 
University of Jyv\"askyl\"a, 40014 Jyv\"askyl\"a, Finland}
\email{bhayo.barkat@gmail.com}

\author[L. Yin]{Li Yin}
\address[L. Yin]{Department of Mathematics, Binzhou University, Binzhou City, Shandong Province, 256603, China}
\email{\href{mailto: L. Yin<yinli_79@163.com>}{yinli\_79@163.com}}

\newcommand{\comment}[1]{}

\swapnumbers
\theoremstyle{plain}

\newtheorem{theorem}[equation]{Theorem}
\newtheorem{theorem a}[equation]{Theorem A}
\newtheorem{theorem b}[equation]{Theorem B}
\newtheorem{lemma}[equation]{Lemma}
\newtheorem{definition}[equation]{Definition}

\newtheorem{remark}[equation]{Remark}

\numberwithin{equation}{section}

\begin{document}

\def\thefootnote{}
\footnotetext{ \texttt{\tiny File:~\jobname .tex,
          printed: \number\year-\number\month-\number\day,
          \thehours.\ifnum\theminutes<10{0}\fi\theminutes}
} \makeatletter\def\thefootnote{\@arabic\c@footnote}\makeatother

\thanks {The second author was supported by the National Natural Science Foundation of China 11401041, and by the
Science Foundation of Binzhou University under grant BZXYL1303}

\begin{abstract}
In this paper, authors study the convexity and concavity properties of real-valued function with respect to the classical means, and prove a conjecture posed by Bruce Ebanks in \cite{e}.
\end{abstract}

\maketitle


{\bf keywords:} Logarithmic mean, identric mean, power mean, Alzer mean, convexity and concavity property, conjecture.

\section{Introduction}

A function $f:\mathbb{R}_+\to\mathbb{R}_+$ is $[m_1,m_2]$-convex (concave) if 
$f(m_1(x,y))\leq (\geq)m_2(f(x),f(y))$ for all $x,y\in\mathbb{R}_+=(0,\infty)$ 
and $m_1,m_2\in \mathbb{M}$, where $\mathbb{M}$ denotes the family of all mean values of two numbers in $\mathbb{R}_+$. Some examples of mean values of two distinct positive real numbers are given below:

\begin{eqnarray*}
\emph{Arithmetic mean:}& &A= A(x,y)=\frac{x+y}{2},\\
\emph{Geometric mean:}& &G= G(x,y)=\sqrt{xy},\\
\emph{Harmonic mean:}& &H=H(x,y)=\frac{1}{A(1/x,1/y)},\\
\emph{Logarithmic mean:}& &L=L(x,y)=\frac{x-y}{\log(x)-\log(y)},\\
\emph{Identric mean:}& &I=I(x,y)=\frac{1}{e}\left( {\frac{{x^x }}{{y^y }}} 
\right)^{1/(x - y)},\\
\emph{Alzer mean:}& &J_p=J_p(x,y)=
\frac{p}{{p + 1}}\frac{{x^{p + 1}  - y^{p + 1} }}{{x^p  - y^p }},\quad p \ne 0, - 1,\\
\emph{Power mean:}& &M_t=M_t(x,y)=\displaystyle\left\{\begin{array}{lll} \displaystyle\left(\frac{x^t+y^t}{2}\right)^{1/t},\;\quad t\neq 0,\\
\sqrt{x\,y}, \,\quad\qquad\qquad t=0\,.\end{array}\right.\\
\end{eqnarray*} 
It is easy to see that $ J_1(x,y)=A(x,y), J_0(x,y)=L(x,y), J_{-2}(x,y)=H(x,y).$ 
For the historical background of these means we refer the reader to see \cite{alzer2,aq, carlson,mitri,sandorc} and the bibliography of these papers.

Before we introduce the earlier results from the literature we recall the following definition, see \cite{avv,aumann}.
\begin{definition}
Let $f:I\to (0,\infty)$ be continuous, where $I$ is a sub-interval of $(0,\infty)$. Let $M$ and $N$ be two any mean functions. We say that the function
$f$ is $MN$-convex (concave) if
$$f (M(x, y)) \leq (\geq)
N(f (x), f (y)) \,\, \text{ for \,\, all} \,\, x,y \in I\,.$$
\end{definition}

In \cite{avv}, Anderson, Vamanamurthy and Vuorinen studied the convexity and concavity of a function $f$ with respect two mean values, and gave the following detailed result:

\begin{theorem} Let $I$ be an open sub-interval of $(0,\infty)$ and let 
$f : I\to (0,\infty)$ be differentiable.
In parts (4)–(9), let $I = (0, b),\, 0<b<\infty$.
\begin{enumerate}
\item $f$ is $AA$-convex (concave) if and only if $f'(x)$ is increasing (decreasing),
\item $f$ is $AG$-convex (concave) if and only if $f'(x)/f (x)$ is increasing (decreasing),
\item $f$ is $AH$-convex (concave) if and only if $f'(x)/f (x)^2$ is increasing (decreasing),
\item $f$ is $GA$-convex (concave) if and only if $xf'(x)$ is increasing (decreasing),
\item $f$ is $GG$-convex (concave) if and only if $xf'(x)/f (x)$ is increasing (decreasing),
\item $f$ is $GH$-convex (concave) if and only if $xf'(x)/f (x)^2$ is increasing (decreasing),
\item $f$ is $HA$-convex (concave) if and only if $x^2f'(x)$ is increasing (decreasing),
\item $f$ is $HG$-convex (concave) if and only if $x^2f'(x)/f (x)$ is increasing (decreasing),
\item $f$ is $HH$-convex (concave) if and only if $x^2f'(x)/f (x)^2$ is increasing 
(decreasing).
\end{enumerate}
\end{theorem}
After the publication of \cite{avv} many authors have studied generalized convexity. For a partial survey of the recent results, see  \cite{avz}.

In \cite{by2}, the following inequalities were studied:

\begin{theorem}\label{3.1-thm}
Let $f:I\to(0,\infty)$ be a continuous and $I\subseteq(0,\infty)$, then 
\begin{enumerate}
\item $f$ is $LL$-convex (concave) if $f$ is increasing and $\log$-convex (concave),
\item $f$ is $AL$-convex (concave) if $f$ is increasing and $\log$-convex (concave),
\end{enumerate}
\end{theorem}

Recently, Baricz \cite{baricz} took one step further and studied the $MN$-convexity(concavity) of a function $f$ in a generalized way, and gave the following result:

\begin{lemma}\cite[Lemma 3]{baricz}\label{bari lemma 3}
  Let $p,q \in \mathbb{R}$ and let $f \colon [a,b] \to (0,\infty)$ be a differentiable function for $a,b \in (0,\infty)$. The function $f$ is $(p,q)$-convex ($(p,q)$-concave) if and only if $x \mapsto x^{1-p} f'(x) (f(x))^{q-1}$ is increasing (decreasing).
\end{lemma}

It can be observed easily that  
$(1, 1)$-convexity means the $AA$-convexity,
$(1, 0)$-convexity means the $AG$-convexity, and $(0, 0)$-convexity means
$GG$-convexity.

\begin{lemma}\cite[Theorem 7]{baricz}\label{bari theorem 7}
  Let $a,b \in (0,\infty)$ and $f \colon [a,b] \to (0,\infty)$ be a differentiable function. Denote $g(x) = \int_1^x f(t) \, dt$ and $h(x) = \int_x^b f(t) \, dt$. Then
  
  \noindent (a) If for all $p \in [0,1]$ the function $f$ is $(p,0)$-concave, then the function $g$ is $(p,q)$-concave for all $p \in [0,1]$ and $q \le 0$. If, in addition the function $x \mapsto x^{1-p} f(x)$ is increasing for all $p \in [0,1]$, then $g$ is $(p,q)$-concave for all $p \in [0,1]$ and $q \in (0,1)$. Moreover, if for all $p \in \mathbb{R}$ the function $x \mapsto x^{1-p} f(x)$ is increasing, then $g$ is $(p,q)$-convex for all $p \in \mathbb{R}$ and $q \ge 1$.
  
  \noindent (b) If for all $p \in [0,1]$ the function $f$ is $(p,0)$-concave, then the function $g$ is $(p,q)$-concave for all $p \in [0,1]$ and $q \le 0$. If, in addition the function $x \mapsto x^{1-p} f(x)$ is decreasing for all $p \in [0,1]$, then $g$ is $(p,q)$-concave for all $p \in [0,1]$ and $q \in (0,1)$. Moreover, if for all $p \in \mathbb{R}$ the function $x \mapsto x^{1-p} f(x)$ is decreasing, then $g$ is $(p,q)$-convex for all $p \in \mathbb{R}$ and $q \ge 1$.
  
  \noindent (c) If for all $p \notin (0,1)$ we have $a^{1-p} f(a) = 0$ and the function $f$ is $(p,0)$-convex, then $g$ is $(p,q)$-convex for all $p \notin (0,1)$ and $q \ge 0$. If, in addition the function $x \mapsto x^{1-p} f(x)$ is increasing for all $p \notin (0,1)$, then $g$ is $(p,q)$-convex for all $p \notin (0,1)$ and $q < 0$.
  
  \noindent (d) If for all $p \notin (0,1)$ we have $b^{1-p} f(b) = 0$ and the function $f$ is $(p,0)$-convex, then $g$ is $(p,q)$-convex for all $p \notin (0,1)$ and $q \ge 0$. If, in addition the function $x \mapsto x^{1-p} f(x)$ is decreasing for all $p \notin (0,1)$, then $g$ is $(p,q)$-convex for all $p \notin (0,1)$ and $q < 0$.
\end{lemma}

In this paper we make a contribution to the subject by giving the following theorems.

\begin{theorem}\label{thm2p1}
Let $f:I \rightarrow (0,\infty)$ and $I\subseteq(0,\infty).$
Then the following inequality holds true:
$$  I(f(x),f(y))\geq f(I(x,y))$$
$$I(f(x),f(y))\leq f(A(x,y)),$$
if the function $f(x)$ is a continuously differentiable, increasing and log-convex(concave).
\end{theorem}

%\begin{corollary}\label{cor2p1}
%Let $b,c>1, z \in\bigl(0,1\bigr)$ and $a,a' \in (0,\infty)$, then
%\begin{equation}I\bigl(F(a,b;c;z), F(a',b;c;z)\bigr) \geq F(I(a,a'),b;c;z).\end{equation}
%\end{corollary}

\begin{theorem}\label{thm2p2}
Let $f$ be a continuous real-valued function on $(0,\infty).$ If $f$ is strictly increasing and convex, then
\begin{equation}
  P_f (x,y)\leq R_f (x,y)
\end{equation}
where $$P_f (x,y) = f\left( {\left( {xy} \right)^{{1 \mathord{\left/
 {\vphantom {1 4}} \right.
 \kern-\nulldelimiterspace} 4}} \left( {\frac{{x + y}}{2}} \right)^{{1 \mathord{\left/
 {\vphantom {1 2}} \right.
 \kern-\nulldelimiterspace} 2}} } \right)
$$
and $$R_f (x,y) = \frac{1}{{y - x}}\int_x^y {f(t)dt}.$$
\end{theorem}

\begin{remark}\rm
In \cite{e}, Ebanks defined
$P_f (x,y)$ and $R_f (x,y),$
and proposed an open problem for a continuous and strictly monotonic real-valued function $f$ on $(0,\infty)$ as follows:

\noindent{Problem}. Does $f$ strictly increasing and convex (or $f''>0$) imply $P_f  \le Rf$ ?
 
It is obvious that the Theorem \ref{thm2p2} gives an affirmative answer to the Ebanks' problem.
\end{remark}

\begin{theorem}\label{2p3}
Let $f:I \rightarrow(0,\infty)$ and $I\subseteq(0,\infty).$ \\
(1)~~If $f(x)$ is continuously differentiable, strictly increasing(decreasing) and convex(concave) and $f^{p-1} (x)f
(x)$ is increasing on $(0,1)$, then\\
$$  J_{p} (f(x),f(y))\geq f(J_{p} (x,y))$$
$$J_{p} (f(x),f(y))\leq f(A(x,y))$$ for $p\leq 1$.\\
(2) ~~If $f(x)$ is continuously differentiable, strictly decreasing(increasing) and convex(concave) and $f^{p-1} (x)f(x)$ is decreasing on $(0,1)$, then\\
$$  J_{p} (f(x),f(y))\geq f(J_{p} (x,y))$$
$$J_{p} (f(x),f(y))\leq f(A(x,y))$$ for $p> 1$.
\end{theorem}

\section{Lemmas and proofs}

We recall the following lemmas which will be used in the proof of the theorems.
\begin{lemma}\cite{qh}\label{2.1-lem}
Let $f,g:[a,b] \to R$ be integrable functions, both increasing or both decreasing. Furthermore, let $p:[a,b] \to R$ be a positive, integrable function. Then
\begin{equation}\label{(2.1)}
\int_a^b {p(x)f(x)dx \cdot } \int_a^b {p(x)g(x)dx}  \le \int_a^b {p(x)dx}  
\cdot \int_a^b {p(x)f(x)g(x)dx}.
\end{equation}
If one of the functions $f$ or $g$ is non-increasing and the other 
non-decreasing, then the inequality in (3.1) is reversed.
\end{lemma}

\begin{lemma}\label{2.2-lem}\cite{kua}
If $f(x)$ is continuous and convex function on $[a,b],$ and $\varphi(x)$ 
is continuous on $[a,b],$ then
\begin{equation}\label{(2.2)}
f\left( {\frac{1}{{b - a}}\int_a^b {\varphi (x)dx} } 
\right) \le \frac{1}{{b - a}}\int_a^b {f\left( {\varphi (x)} \right)dx}.
\end{equation}
If function $f(x)$ is continuous and concave on $[a,b],$ 
the inequality in \eqref{(2.2)} is reversed.
\end{lemma}

\begin{lemma}\cite{aq}\label{2.3-lem}
Fix two positive number $a, b.$
Then $L(a,b)\leq I(a,b)\leq A(a,b).$
\end{lemma}

\begin{lemma}\label{lem34}\cite{kua}
The function $p\mapsto J_p(x,y)$ is strictly increasing on 
$\mathbb{R}\setminus\{0,-1\}$.
\end{lemma}

%\begin{lemma}\cite[Theorem 2]{bbv}
%For all $x \in(0,\infty )$ fixed, The following hold\\
%(4)~~ $p\mapsto {\rm arcsinh}_{p,q}(x)$ is strictly increasing and concave on $(1,\infty)$ for $q>1$.\\
%(5)~~ $q\mapsto {\rm arcsinh}_{p,q}(x)$ is strictly increasing and concave on $(1,\infty)$ for $p>1$.\\
%\end{lemma}
%
%\begin{lemma}\cite[Theorem 1]{bbk}
%If $p,q>1$ and $a \geq1$, then ${\rm arcsin}_{p,q}(x)$ is $(a,a)-$convex on $(0,1)$, ${\rm arctan}_{p,q}(x)$ is $(a,a)-$concave on $(0,1)$, while ${\rm arcsinh}_{p,q}(x)$ is $(a,a)-$concave on $(0,\infty)$. \end{lemma}
%
%\begin{lemma}\cite[Theorem 2]{bbk}
%If $p,q>1$ and $a \geq1$, then ${\rm sin}_{p,q}(x)$ is $(a,a)-$concave on $(0,1)$, and ${\rm cos}_{p,q}(x), {\rm tan}_{p,q}(x), {\rm sinh}_{p,q}(x)$ is $(a,a)-$ convex on $(0,1)$. \end{lemma}

%\begin{lemma}\cite{neu}
%For $b,c>0$ and $x \in(0,1)$, the function 
%$a\mapsto F\left(a,b;c;x \right)$ is increasing and 
%logarithmically convex on $(0,\infty)$. 
%\end{lemma}

%\section{Proofs of the main results}
\noindent{\bf Proof of Theorem \ref{thm2p1}.}
Since the proof of part (2) is similar to part (1), 
we only prove the part (1) here.
An easy computation and substitution $t=f(u)$ yield
\begin{equation}\label{eq1a}
\ln I(f(x),f(y)) = \frac{{\int_{f(y)}^{f(x)} {\ln tdt} }}
{{\int_{f(y)}^{f(x)} {1} }} = \frac{{\int_y^x { \ln f(u)f'(u)du} }}
{{\int_y^x { f'(u) du} }}.
\end{equation}
Since the functions $f(x)$ and  $ f'(x) $ are increasing 
on $I\subseteq(0,\infty)$, now by using Lemma \ref{2.1-lem} we have
\begin{equation}\label{eq2a}
 \int_y^x {1du}  \cdot \int_y^x {\ln f(u)f'(u)du}  
\geq \int_y^x { f'(u) du}  
\cdot \int_y^x {\ln f(u)du}.
\end{equation}
Combining \eqref{eq1a} and \eqref{eq2a},
we obtain
$$
I(f(x),f(y)) \geq \frac{{\int_y^x {\ln f(u)du} }}{{y - x}}.
$$
Considering the log-convexity of the function $f(x)$ and using Lemmas \ref{2.2-lem} and 
\ref{2.3-lem}, we get
$$
I(f(x),f(y)) \geq \ln f\left( {\frac{{\int_y^x {udu} }}{{y - x}}} \right) 
= \ln f\left( {\frac{{x + y}}{2}} \right) \geq \ln f\left( {I(x,y)} \right).
$$
This completes the proof.$\hfill\square$

%\begin{remark}
%Using Lemma 3.5, 3.6, 3.7 3.8 and Theorem 2.1, we easily obtain Corollary 2.1, 2.2, 2.3.
%\end{remark}
\vspace{.3cm}
\noindent{\bf Proof of Theorem \ref{thm2p2}.}
Since $f$ is strictly increasing and convex, by utilizing the 
Lemma \ref{2.1-lem} and the inequality
$G(x,y)\leq A(x,y)$ we obtain
\begin{eqnarray*}
 R_f (x,y) &\ge& \frac{{\int_x^y {f(u)du} }}{{y - x}} \ge f\left( {\frac{{\int_x^y {udu} }}{{y - x}}} \right) \\
  &=& f\left( {\frac{{x + y}}{2}} \right) \ge f\left( {\left( {xy} \right)^{1/4} \left( {\frac{{x + y}}{2}} \right)^{1/2} } \right)\\
	&=& P_f (x,y).
\end{eqnarray*}
This completes the proof. $\hfill\square$

%\begin{remark}
%Theorem \ref{thm2p2} solves an open problem proposed by B. Ebanks in \cite{e}.
%\LARGE{Please give details here, like what is the problem and how it can be solved}
%\end{remark}

%\begin{remark}
%Using Lemma 3.6, 3.7 and Theorem 2.2, we easily obtain Corollary 2.4.
%\end{remark}
\vspace{.3cm}
\noindent{\bf Proof of Theorem \ref{2p3}.} For the proof of part (1), 
letting $t=f(u)$, we get
$$
J_{p} (f(x),f(y)) = \frac{{\int_{f(y)}^{f(x)} {t^p dt} }}{{\int_{f(y)}^{f(x)} {t^{p-1}} }} = \frac{{\int_y^x { f^{p} (u)f'(u)du} }}{{\int_y^x { f^{p-1} (u)f'(u) du} }}.
$$
By using Lemma \ref{2.1-lem},
we obtain
$$
J_{p} (f(x),f(y)) \geq \frac{{\int_y^x { f(u)du} }}{{y - x}}.
$$
Considering convexity of the function $f(x)$ and using Lemmas \ref{2.2-lem} 
and \ref{lem34}, we get
$$
J_{p} (f(x),f(y)) \geq f\left( {\frac{{\int_y^x {udu} }}{{y - x}}} \right) 
= f\left( {\frac{{x + y}}{2}} \right) \geq f\left( {J_{p} (x,y)} \right),
$$
this implies (1). The proof of part (2) follows similarly.
$\hfill\square$

\end{document}